\documentstyle[epsf,amsfonts,amssymb]{article} 

\def\abstract#1{\vskip 7mm 
        \begin{center}{\large Abstract}\par \smallskip
                \begin{minipage}[c]{12cm}
                        \small #1
                \end{minipage}
        \end{center}
}
\def\title#1{\begin{center}{\Large\bf #1}\end{center}}
\def\author#1{\vskip 5mm \begin{center}{#1}\end{center}}
\def\address#1{\begin{center}{\it #1}\end{center}}

\def\ben{\begin{equation}}
\def\een{\end{equation}}
\def\bena{\begin{eqnarray}}
\def\eena{\end{eqnarray}}

\newtheorem{proposition} {Proposition} 
\newtheorem{definition}  {Definition}

\newtheorem{corollary}   {Corollary}

\makeatletter
\@ifundefined{lesssim}{}{}
\@ifundefined{gtrsim}{}{}
\def\vereq#1#2{\lower3pt\vbox{\baselineskip1.5pt \lineskip1.5pt
\ialign{$\m@th#1\hfill##\hfil$\crcr#2\crcr\sim\crcr}}}
\makeatother

\begin{document}


\title{%
   Convex Functions and Spacetime Geometry\footnote{
First appeared as pp. 227-230
of Proceedings of the International Workshop: Frontiers of Cosmology and
Gravitation held at the Yukawa Institute, April 25-27, 2001 ed. M. Sakagami.
}
}

\author{%
  Gary W. Gibbons,\footnote{The current e-mail:gwg1@damtp.cam.ac.uk }
}
\address{%
DAMTP, Center for Mathematical Sciences,
University of Cambridge, \\
Wilberforce Road, Cambridge CB3 0WA, the United Kingdom 
}
\author{%
  Akihiro Ishibashi\footnote{
The current e-mail and affiliation: akihiro@phys.kindai.ac.jp \, 
Kindai University, Higashi-Osaka, Japan 
} 
}
\address{
  Yukawa Institute for Theoretical Physics, Kyoto University, \\ 
  Sakyo-ku, Kyoto 606--8502, Japan
}

\abstract{
Convexity and convex functions play an important role
in theoretical physics. 
To initiate a study of the possible uses of convex functions 
in General Relativity, 
we discuss the consequences of a spacetime $(M,g_{\mu \nu})$ or 
an initial data set $(\Sigma, h_{ij}, K_{ij})$ admitting 
a suitably defined convex function. 
We show how the existence of a convex function on a spacetime 
places restrictions on the properties of the spacetime geometry. 
}

\section*{Summary}

Convexity and convex functions play an important role
in theoretical physics. For example, Gibbs's approach to 
thermodynamics~\cite{Gibbs} 
is based on the idea that the free energy 
should be a convex function. 
A closely related concept is that of a convex cone
which also has numerous applications to physics. 
Perhaps the most familiar example is the lightcone of Minkowski spacetime.
Equally important is the convex cone of mixed states
of density matrices in quantum mechanics. 
Convexity and convex functions also have important applications to
geometry, including Riemannian geometry~\cite{Udriste.C1994}. 
It is surprising therefore that, to our knowledge, that techniques making 
use of convexity and convex functions have played no great role 
in General Relativity. 
The purpose of this paper is to initiate a study of the possible
uses of such techniques. 

We give the definition of a convex function on a spacetime $(M, g_{\mu \nu})$ 
as follows.

\begin{definition}{\bf Spacetime definition:} \\
A smooth function $f: M \rightarrow {\Bbb R}$ is called a 
{\em spacetime convex function} if the Hessian $\nabla_\mu \nabla_\nu f$ has 
Lorentzian signature and satisfies the condition, 
\ben
V^\mu V^\nu \nabla_\mu \nabla_\nu f \ge c g_{\mu \nu} V^\mu V^\nu \;, \,\,\, 
for \,\, {}^\forall V^\mu \in TM \;, \,\, {} c>0 \, \, constant. 
\label{def:spacetimeconvex}
\een 
\end{definition} 

\noindent
Since a spacetime convex function has a Hessian with Lorentzian signature, 
we can say that the spacetimes admitting spacetime convex functions 
have particular types of causal structures. 

One of the simplest examples of a spacetime convex function 
is the canonical one, 
\ben
  f = \frac{1}{2} \left( x^i x^i - \alpha t^2 \right) \;, 
          \quad (t,x^i) \in M \;, 
\label{convex:canonical}
\een
where $\alpha$ is a constant such that $0 < \alpha \le 1$. 

\bigskip 

The existence of a convex function, suitably defined, on a spacetime 
places important restrictions on the properties of the spacetime. 
For example, we have the following propositions. 

\bigskip 

Suppose a spacetime convex function $f$ exists in $(M,g_{\mu \nu})$. 

\begin{proposition}
 $(M,g_{\mu \nu})$ admits no closed spacelike geodesics. 
\end{proposition}

\begin{proposition} 
Consider a spacetime $(M,g_{\mu \nu})$ and a closed spacelike surface $S \subset M$. 
Suppose that in a neighbourhood of $S$ the metric is written as 
\ben 
g_{\mu \nu}dx^\mu dx^\nu = \gamma_{ab} dy^a dy^b + k_{pq}dz^p dz^q \;, 
\een 
where $k_{pq}dz^p dz^q$ is the metric on $S$ and 
the components of the two-dimensional Lorentzian metric $\gamma_{ab}$ are 
independent of the coordinates $z^p$. 
Then, $S$ cannot be a closed marginally inner and outer trapped surface.  
\label{proposition:trappedsurface}
\end{proposition}

Assuming that $(M,g_{\mu \nu})$ admits a spacetime convex function $f$, 
we deduce 
\begin{proposition} 
 If $\Sigma$ is totally geodesic, i.e., $\Sigma$ is a surface of time 
symmetry, then $(\Sigma, h_{ij})$ admits a convex function.  
\end{proposition}

\begin{proposition} 
 If $\Sigma$ is maximal, i.e., $h^{ij}K_{ij} = 0$, it admits 
 subharmonic function. 
\end{proposition}
\begin{corollary}
 $(\Sigma, h_{ij})$ cannot be closed. 
\label{coro:cannot be closed}
\end{corollary}
\begin{corollary}
 $(\Sigma, h_{ij})$ admits no closed geodesics nor minimal two surface 
in the case that $K_{ij}=0$. 
\label{corollary:minimalsurface}
\end{corollary} 

\bigskip 

One way of viewing the ideas of this paper is in terms of a sort of duality
between paths and particles on the one hand and functions and waves 
on the other. Mathematically the duality corresponds
to interchanging range and domain. A curve $x(\lambda)$ is a map 
$x: {\Bbb R} \rightarrow M$ while a function $f(x)$ 
is a map $f:M \rightarrow {\Bbb R}$. 
A path arises by considering invariance under
diffeomorphisms of the domain (i.e. of the world volume)
and special paths, for example geodesics have action functionals
which are reparametrization invariant. Much effort, physical and
mathematical has been expended on exploring the global
properties of spacetimes using geodesics. 
Indeed there is a natural
notion of convexity based on geodesics. 
Often a congruence of geodesics is used.

On the dual side, one might consider properties
which are invariant  under diffeomorphisms $f(x)\rightarrow g(f(x))$
of the range or target space. That is one 
may explore  the global properties
of spacetime using the foliations provided
by the level sets of a suitable function. The analogue
of the action functional for geodesics is one like
\ben
\int_M \sqrt {\epsilon \nabla _\mu f \nabla^\mu f}, 
\een
where $\epsilon =\pm$, depending upon whether $\nabla _\mu  f$ is 
spacelike or timelike, and which is invariant under reparametrizations 
of the range. 
The level sets  of the solutions of the Euler-Lagrange equations
are then minimal surfaces. One may also consider
foliations by totally umbilic surfaces or by ``trace $K$ equal constant" 
foliations, and this is often done in numerical relativity.
The case of a convex function then corresponds 
to a foliation by totally expanding hypersurfaces, 
that is, hypersurfaces with positive definite second fundamental form.

Actually, given two vectors $X^\mu$ and $Y^\nu$ tangent to a level set 
$\Sigma_c=\{ x \in M | f(x) = constant = c\}$ of a function $f$, 
one may evaluate $K_{\mu \nu}X^\mu Y^\mu$ in terms of the Hessian of $f$:
\ben
K_{\mu \nu}X^\mu Y^\mu= X^\mu Y^\nu {\nabla _\mu \nabla _\nu 
f \over \sqrt{\epsilon \nabla _\nu f
\nabla ^\nu f}}.
\een
Thus, if $M$ is a Riemannian manifold, a strictly convex function has a positive definite 
second fundamental form. 
Of course there is a convention here about the choice
of direction of the normal. We have chosen $f$ to decrease along
$n^\mu$. The converse is not necessarily true,
since $f(x)$ and $g(f(x))$ have the same level sets, where $g$ is 
a monotonic function of ${\Bbb R}$. Using this gauge freedom 
we may easily change the signature  of the Hessian.
However, given a hypersurface $\Sigma_0$ 
with positive definite second fundamental form,
we may, locally,  also use this gauge freedom to find 
a convex function  whose level $f=0$  coincides with $\Sigma_0$. 
If we have a foliation (often called a ``slicing"   
by relativists) by hypersurfaces with positive definite fundamental
form, we may locally represent the leaves as the levels sets of a 
convex function. 
   
Similar remarks apply for  Lorentzian metrics.
The case of greatest interest is when the level sets have a timelike
normal.  For a classical strictly convex function, the second fundamental 
form will be positive definite and the hypersurface orthogonal timelike
congruence given by the normals $n^\mu$ is an expanding one.
For a spacetime convex function, our conventions also imply that
the second fundamental form of a spacelike level set 
is positive definite. This can be illustrated by the canonical
example (\ref{convex:canonical}) with $\alpha =1$ in flat spacetime. 
The spacelike level sets foliate the interior of the future (or past) 
light cone. Each leaf is isometric to hyperbolic space. 
The expansion is homogeneous and isotropic, and
if we introduce coordinates adapted to the foliation we
obtain the flat metric in $K=-1$ FLRW form with scale factor $a(\tau)=\tau$. 
This is often called the Milne model.
From what has been said it is clear that there is 
a close relation between the existence of convex functions
and the existence of foliations with positive definite second 
fundamental form. 

\bigskip 

In the theory of maximal hypersurfaces and more generally
hypersurfaces of constant mean curvature (``${\rm Tr}K=constant$" hypersurface) 
an important role is played by the idea of a ``barriers.'' 
The basic idea \cite{Eschenburg1989} is that, if two spacelike hypersurfaces 
$\Sigma_1$ and $\Sigma_2$ touch, then the one in the future, 
$\Sigma_2$ say, can have no smaller a mean curvature than the 
$\Sigma_1$, the hypersurface in the past, i.e., 
${\rm Tr} K_1 \le {\rm Tr} K_2$. 

As an application of this idea, consider the interior region of the 
Schwarzschild solution. 
The hypersurface $r=$ constant has 
\ben
   {\rm Tr} K(r) = - {2 \over r}\left({2M \over r} -1 \right)^{-1/2} 
                          \left(1 - {3M \over 2r} \right) \;. 
\een 
The sign of ${\rm Tr} K(r)$ is determined by the fact that regarding $r$ as
a time coordinate, $r$ decrease as time increases. Thus for $ r >
3M/2$, ${\rm Tr} K$ is negative but for $r <3M/2$, ${\rm Tr} K$ is positive. 
If $r = 3M/2$, we have ${\rm Tr} K=0$, that is, $r = 3M/2$ is a maximal
hypersurface. 

Consider now attempting to foliate the black hole interior region by the level
sets of a concave function. Every level set must touch 
an $r =$ constant hypersurface at some point. 
If $r = 3M/2$ this can, by \cite{Eschenburg1989}, only happen if the level set 
lies in the past of the hypersurface $r = 3M/2$. 
Evidently therefore a foliation by level sets of a spacetime 
{\it concave} function can never penetrate the barrier at $r = 3M/2$. 
In particular, such a foliation can never extend to the singularity at
$r=0$. 

These results are relevant to work in numerical relativity. One
typically sets up a coordinate system in which the constant time
surfaces are of constant mean curvature. 
It follows from our results that, if the constant is {\sl negative}, 
then, the coordinate system can never penetrate the region $r<3M/2$. 
In fact it could never penetrate any maximal hypersurface.  
However, level sets $\Sigma$s of a convex function can have positive 
${\rm Tr} K$, hence can penetrate the barrier surface.

\bigskip 

The existence problem of constant mean curvature foliations 
has been investigated extensively not only in black hole spacetimes 
as discussed here but also in cosmological spacetimes. 
In cosmological spacetimes, constant mean curvature hypersurfaces, if exist, 
are likely to be compact, and thus do not admit a strictly or uniformly 
convex function which live on the hypersurfaces, because of 
Corollary~\ref{coro:cannot be closed}. 
However, a {\sl spacetime convex} function, if available, can give a 
constant mean curvature foliation with non-vanishing mean curvature 
as its level surfaces. 

\bigskip 

We anticipate that study of convex functions and foliations 
by convex surfaces will provide further insights into global problems 
in general relativity and should have applications to numerical relativity. 

\bigskip 

We have mainly discussed the relations between submanifolds such as
hypersurfaces in a spacetime and the existence of convex functions 
in the ambient manifold. 
The proofs of Propositions in this article, and further examples of 
spacetime convex functions are given in the paper~\cite{GI}. 

\newpage 

\section*{\bf Acknowledgments}
G.W.G. would like to thank Profs. T. Nakamura and H. Kodama for their 
kind hospitality at the Yukawa Institute for Theoretical Physics 
where the main part of this work was done. 
A.I. would like to thank members of DAMTP for their kind hospitality. 
A.I. was supported by Japan Society for the Promotion of Science.  


\end{document}